\documentclass[12pt]{article}

\usepackage{amsmath,amssymb}

\begin{document}

\pagestyle{myheadings} \markright{LEFSCHETZ FORMULAE...}

\def \1{{\bf 1}}
\def \a{{{\mathfrak a}}}
\def \ad{{\rm ad}}
\def \al{\alpha}
\def \ar{{\alpha_r}}
\def \A{{\mathbb A}}
\def \Ad{{\rm Ad}}
\def \Aut{{\rm Aut}}
\def \b{{{\mathfrak b}}}
\def \bs{\backslash}
\def \B{{\cal B}}
\def \c{{\mathfrak c}}
\def \cent{{\rm cent}}
\def \C{{\mathbb C}}
\def \CA{{\cal A}}
\def \CB{{\cal B}}
\def \CC{{\cal C}}
\def \CD{{\cal D}}
\def \CE{{\cal E}}
\def \CF{{\cal F}}
\def \CG{{\cal G}}
\def \CH{{\cal H}}
\def \CHC{{\cal HC}}
\def \CL{{\cal L}}
\def \CM{{\cal M}}
\def \CN{{\cal N}}
\def \CP{{\cal P}}
\def \CQ{{\cal Q}}
\def \CO{{\cal O}}
\def \CS{{\cal S}}
\def \CT{{\cal T}}
\def \CV{{\cal V}}
\def \d{{\mathfrak d}}
\def \det{{\rm det}}
\def \df{\ \begin{array}{c} _{\rm def}\\ ^{\displaystyle =}\end{array}\ }
\def \diag{{\rm diag}}
\def \dist{{\rm dist}}
\def \End{{\rm End}}
\def \eps{\varepsilon}
\def \eqn{\begin{eqnarray*}}
\def \endeqn{\end{eqnarray*}}
\def \F{{\mathbb F}}
\def \Fix{{\rm Fix}}
\def \Frob{{\rm Frob}}
\def \Fx{{\mathfrak x}}
\def \FX{{\mathfrak X}}
\def \g{{{\mathfrak g}}}
\def \ga{\gamma}
\def \Ga{\Gamma}
\def \Gal{{\rm Gal}}
\def \GL{{\rm GL}}
\def \h{{{\mathfrak h}}}
\def \Hom{{\rm Hom}}
\def \im{{\rm im}}
\def \ind{{\rm ind}}
\def \Im{{\rm Im}}
\def \Ind{{\rm Ind}}
\def \k{{{\mathfrak k}}}
\def \K{{\cal K}}
\def \l{{\mathfrak l}}
\def \la{\lambda}
\def \lap{\triangle}
\def \li{{\rm li}}
\def \La{\Lambda}
\def \m{{{\mathfrak m}}}
\def \mod{{\rm mod}}
\def \n{{{\mathfrak n}}}
\def \name{\bf}
\def \Mat{{\rm Mat}}
\def \N{\mathbb N}
\def \o{{\mathfrak o}}
\def \ord{{\rm ord}}
\def \O{{\cal O}}
\def \p{{{\mathfrak p}}}
\def \ph{\varphi}
\def \prf{\noindent{\bf Proof: }}
\def \per{{\rm per}}
\def \Per{{\rm Per}}
\def \q{{\mathfrak q}}
\def \qed{\ifmmode\eqno $\square$\else\noproof\vskip 12pt plus 3pt minus 9pt \fi}
 \def\noproof{{\unskip\nobreak\hfill\penalty50\hskip2em\hbox{}
     \nobreak\hfill $\square$\parfillskip=0pt
     \finalhyphendemerits=0\par}}
\def \Q{\mathbb Q}
\def \res{{\rm res}}
\def \R{{\mathbb R}}
\def \Re{{\rm Re \hspace{1pt}}}
\def \r{{\mathfrak r}}
\def \ra{\rightarrow}
\def \rank{{\rm rank}}
\def \span{{\rm span}}
\def \supp{{\rm supp}}
\def \SL{{\rm SL}}
\def \Spin{{\rm Spin}}
\def \t{{{\mathfrak t}}}
\def \T{{\mathbb T}}
\def \tr{{\hspace{1pt}\rm tr\hspace{2pt}}}
\def \vol{{\rm vol}}
\def \z{\zeta}
\def \Z{\mathbb Z}
\def \={\ =\ }

\newcommand{\frack}[2]{\genfrac{}{}{0pt}{}{#1}{#2}}
\newcommand{\rez}[1]{\frac{1}{#1}}
\newcommand{\der}[1]{\frac{\partial}{\partial #1}}
\renewcommand{\binom}[2]{\left( \begin{array}{c}#1\\#2\end{array}\right)}
\newcommand{\norm}[1]{\left|\hspace{-1pt}\left| #1\right|\hspace{-1pt}\right|}
\renewcommand{\matrix}[4]{\left(\begin{array}{cc}#1 & #2 \\ #3 & #4 \end{array}\right)}
\renewcommand{\sp}[2]{\langle #1,#2\rangle}

\newtheorem{theorem}{Theorem}[section]
\newtheorem{conjecture}[theorem]{Conjecture}
\newtheorem{lemma}[theorem]{Lemma}
\newtheorem{corollary}[theorem]{Corollary}
\newtheorem{proposition}[theorem]{Proposition}

\title{Lefschetz formulae and zeta functions}
\author{Anton Deitmar}

\date{}
\maketitle

$$ $$

{\bf Abstract:}
The connection between Lefschetz formulae and zeta function is explained. As a particular example the theory of the generalized Selberg zeta function is presented. Applications are given to the theory of Anosov flows and prime geodesic theorems.



\tableofcontents

$$ $$

\begin{center} {\bf Introduction} \end{center}

This is a survey on the connection between Lefschetz formulae and zeta functions in
general and the work of the author on the generalized Selberg zeta function in particular.

A Lefschetz formula relates fixed points of an automorphism $f$ of some space to some global
cohomology and the action of $f$ thereon.
In the realm of geometrically defined zeta functions they serve to prove rationality by giving an
interpretation of a zeta function as a ``characteristic function'' of the induced automorphism on the
cohomology.

Such an interpretation is highly desirable for other types of zeta functions such as the Riemann zeta
function. Since the Riemann zeta function is not rational the cohomology space in question should
be infinite dimensional. In this paper we give some prototypes of zeta functions which can be
interpreted in this way for some natural infinite dimensional cohomology groups.

Starting with discrete dynamical systems for which the classical Lefschetz formula ensures
rationality of the zeta function, the case of suspensions, which is
treated in some detail, serves as a guideline to find the suitable cohomology theory. It turns out that
foliation cohomology, or rather reduced foliation cohomology does the trick. 

The Selberg zeta function is an example of a zeta function which can be interpreted in this way. The
corresponding Lefschetz formula can be generalized to a multi-dimensional Lefschetz formula for
higher rank spaces. The latter can be applied to give meromorphic extension of generalized Selberg
zeta functions or to prove higher rank prime geodesic theorems.

$ $

\section{The classical Lefschetz formula}

Let $M$ be a compact smooth manifold and let $f:M\ra M$ be a diffeomorphism. We say that $f$ is
\emph{regular} if its graph $\Ga(f)$ intersects the diagonal $\Delta\subset M\times M$ transversally
only. This is equivalent to saying that for every fixed point $x$ of $f$ the differential $f_*:T_xM\ra
T_xM$ satisfies $\det(1-f_*\mid T_xM)\ne 0$. In that case we define the \emph{index} of the fixed
point $x$ as
$$
\ind_f(x)\df {\rm sign}\left(\det(1-f_*\mid T_xM)\right)\ \in\ \{\pm 1\}.
$$
Note that the regularity implies that the fixed points of $f$ are isolated. Since $M$ is compact they
are finite in number. 

\begin{theorem}\label{Lefschetz classical} (Lefschetz trace formula)\\
If $f$ is regular, then
$$
\sum_{x=f(x)}\ind_f(x)\=\sum_{q=0}^{\dim M} (-1)^q\,\tr(f^*\mid H^q(M)).
$$
Here $H^q(M)$ denotes the cohomology of $M$ with, say, complex coefficients.
\end{theorem}

For a proof see \cite{Dieu}.

Suppose that every iterate $f^n=f\circ\cdots\circ f$ of $f$ has only finitely many fixed points. In that
case define the zeta function of $f$ as
$$
Z_f(T)\=\exp\left(\sum_{n=1}^\infty \frac{T^n}{n}\,\# \Fix(f^n)\right),
$$ 
where $\Fix(f^n)$ is the set of fixed points of the map $f^n$. 
Since $\# \Fix(f^n)=\sum_{d|n}\sum_{|o|=d} d$, where the outer sum runs over all positive
divisors of $n$ and the inner sum runs over all $f$-orbits $o$ in $M$ of cardinality $d$. This implies that
$$
Z_f(T)\=\prod_{o}(1-T^{|o|})^{-1},
$$
where the product runs over all finite $f$-orbits in $M$. At first this is a formal power series in $T$.
Under additional assumptions we can say more.

We say that $f$ is \emph{strongly regular} if all iterates $f^n$, $n\in\N$ are regular and if
$\ind_{f^n}(x)=1$ for every $x\in \Fix(f^n)$. The latter condition follows for example if $M$ is a complex
manifold and $f$ is holomorphic.

\begin{theorem}(Lefschetz determinant formula)\\
If $f$ is strongly regular, then
$$
Z_f(T)\=\prod_{q=0}^{\dim M}\det\left( 1-Tf^*\mid H^q(M)\right)^{(-1)^{q+1}}.
$$
In particular, $Z_f(T)$ is a rational function in $T$.
\end{theorem}

\prf This is a simple consequence of Theorem \ref{Lefschetz classical}.
The assumptions imply that $\#\Fix(f^n)=\sum_{x=f^n(x)}\ind_{f^n}(x)$, so the Lefschetz formula gives
\begin{eqnarray*}
Z_f(T) &=& \exp\left(\sum_{n=1}^\infty \frac{T^n}n \sum_{q=0}^{\dim M} (-1)^q\,\tr(f^* \mid H^q(M))\right)\\
&=& \prod_{q=0}^{\dim M} \exp\left((-1)^q\,\tr\sum_{n=1}^\infty \frac{(Tf^*)^n}{n} \mid H^q(M)\right)\\
&=& \prod_{q=0}^{\dim M}\det\left(\exp(-\log(1-Tf^*)\mid H^q(M)\right)^{(-1)^{q+1}}\\
&=& \prod_{q=0}^{\dim M} \det(1-Tf^*\mid H^q(M))^{(-1)^{q+1}}.
\end{eqnarray*}
\qed

Lefschetz formulae of the trace or the determinant type emerge in varying contexts throughout
mathematics. As an example consider the Hasse-Weil zeta function $Z_V(T)$ of a smooth projective
variety $V$ over a finite field $\F$. It is defined as
$$
Z_V(T)\=\exp\left(\sum_{n=1}^\infty \frac{T^n}{n} \# V(\F_{n})\right),
$$
where $\F_n$ is the extension of $\F$ of degree $n$, which is uniquely determined up to isomorphism. For this zeta function M.Artin, Verdier, and
Grothendieck showed that
$$
Z_V(T)\=\prod_{j=0}^{2\dim M}\det(1-T\, \Frob\mid H_{l-ad}^j(V))^{(-1)^{j+1}},
$$
where $\Frob$ is the Frobenius acting on the $l$-adic cohomology $H_{l-ad}^\bullet(V)$ for a prime $l$
different from the characteristic of $\F$.

\section{The classical Lefschetz formula for suspensions}
Let $Y=(M\times \R)/\Z$, where $\Z$ acts on $M\times \R$ by $k.(m,s)=(f^k(m),s-k)$. Then $Y$ is called
the \emph{suspension} of $M$. There is a natural flow $\phi_t$ on $Y$ given by
$$
\phi_t[m,s]\= [m,s+t].
$$
This is the suspended flow of $f$. Since any closed orbit of $\phi_t$ must sit over a periodic point of $f$, the closed orbits of $\phi$ have integer lengths and for
$n\in\N$ the closed orbits of length $n$ correspond to $f$-orbits of fixed points of $f^n$. For a closed orbit $c$ let $l(c)$ denote the length of the orbit. Note here that we use the differential geometric notion of closed orbits rather than the group theoretical one. To make this precise define a \emph{periodic point} to be an element $(y,t)$ of $Y\times 90,\infty)$ such that $\phi_t(y)=y$. The period $t$ is not necessarily the prime period of $y$, ie, the smallest $t>0$ with $\phi_t(y)=y$. A \emph{closed orbit} is an equivalence class of periodic points where two periodic points $(y,t)$ and $(z,s)$ are called equivalent if $s=t$ and $\phi_r(y)=z$ for some $r\in\R$. A closed orbit $(y,t)$ is called \emph{primitive} if $t$ is the prime period of $y$.

The projection to the second variable
$$
\pi : Y\ra \R/\Z
$$
is a fibre bundle with fibre $M$. Let $T_M\subset TY$ be the sub-bundle of vectors tangent to fibres.
These vectors are also called the \emph{vertical} vectors. For
$p\ge 0$ let $H^p(T_M)$ be the cohomology of the complex $\Ga^\infty(\bigwedge^\bullet T_M^*)$ with
the natural exterior differential. Then $H^p(T_M)$ is the space of sections of a vector bundle $E^p$
over $\R/\Z$ whose fibre is $H^p(M)$. This vector bundle can be described as
$$
E^p\= (H^p(M)\times\R)/\Z,
$$
where $\Z$ acts by
$$
k.(v,x)\=(f^{*k}v, x-k).
$$
A section $s$ of $E^p$ can be viewed as a map $s:\R\ra H^p(M)$ satisfying $s(x+k)=f^{*k}s(x)$. The
flow acts on $s$ by $\phi_t^* s(x)\=s(x+t)$.

Let $\ph$ be a smooth function of compact support in $(0,\infty)$. 
Consider the operator
$L_\ph^p=\int_0^\infty \ph(t)\, \phi_t^*\, dt$ on the space of $L^2$-sections of $E^p$.

\begin{theorem} (Lefschetz Formula for suspensions)\\
$L_\ph^p$ is of trace class and
$$
\sum_{p=0}^{\dim M} (-1)^p\, \tr L_\ph^p\=\sum_{c\ \rm closed} l(c_0)\, \ind(c)\, \ph(l(c)),
$$
where the sum on the right hand side runs over all closed orbits of $\phi_t$ and $c_0$ is the primitive
closed orbit underlying $c$.\\
 Further, $\ind(c)={\rm sign}(\det(1-\phi_{l(c), *}\mid T_{M})$.
\end{theorem}

The theorem can be reformulated as an identity of distributions on $(0,\infty)$ as follows.

\begin{corollary}
As an identity of distributions on $(0,\infty)$ we have
$$
\sum_{p=0}^{\dim M} (-1)^p\, \tr (\phi_t^*\mid H^p(T_M))\=\sum_{c\ \rm closed} l(c_0)\, \ind(c)\,
\delta(t-l(c)),
$$
where $\delta$ is the delta distribution.
\end{corollary}

\prf (of the Theorem)
Let $s$ be a section of $E^p$. Then
\begin{eqnarray*}
L_\ph^p\, s(x) &=& \int_\R \ph(t)\, s(x+t)\, dt\\
&=& \int_\R \ph(t-x)\, s(t)\, dt\\
&=&\sum_{k\in\Z} \int_0^1\ph(t-x+k)\, s(t+k)\, dt\\
&=&\sum_{k\in\Z} \int_0^1\ph(t-x+k)\, f^{*k}s(t)\, dt\\
&=& \int_{\R/\Z} K(x,t)\, s(t)\, dt,
\end{eqnarray*}
where $K(x,t)$ is the kernel
$$
K(x,t)\=\sum_{k\in\Z} \ph(t-x+k)f^{*k}.
$$
This sum is locally finite and so $K(x,t)$ is a smooth kernel.
It follows that the operator $L_\ph^p$ is of trace class and that
\begin{eqnarray*}
\tr\, L_\ph^p &=& \int_0^1\tr K(x,x)\, dx\\
&=& \sum_{k\in\Z}\ph(k)\, \tr f^{*k}.
\end{eqnarray*}
So that
\begin{eqnarray*}
\sum_{p=0}^{\dim M}(-1)^p\, \tr L_\ph^p &=& \sum_{k\in\Z} \ph(k)\sum_{p=0}^{\dim M} (-1)^p\,
\tr(f^{*k}\mid H^p(M))\\
&=& \sum_{k\in\N} \ph(k)\sum_{x=f^k(x)}\ind_f(x).
\end{eqnarray*}
Since any closed orbit $c$ of $\phi_t$ gives $l(c_0)$ points $x$ with $x=f^k(x)$ with $k=l(c)$, the claim
follows.
\qed

\section{Foliation cohomology}
The receptacle for the global side of the Lefschetz formula will be a foliation cohomology, a term to be
defined in this section. See also \cite{MooreSchochet}.  A smooth \emph{foliation} on a manifold $M$ is a smooth atlas consisting of
coordinates $(x,y)$ with values in $\R^k\times \R^l$. A set of the form $\{ y\equiv {\rm const}\}$ is
called a \emph{patch} of the coordinate chart. The defining property of a foliation is that 
coordinate changes within the atlas 
map patches to patches. Thus a patch continues into a neighbouring coordinate set and thus
extends to a
$k$-dimensional immersed sub-manifold, called a \emph{leaf} of the foliation. As an example
consider $M=\R^2/\Z^2$. Fix
$\al,\beta\in\R^\times$ and consider the foliation on $M$ with leaves
$$
L_{x,y}\= (x,y)+\R(\al,\beta)\ \ \ \mod\,\Z^2,
$$
then $L_{x,y}$ is the leaf through the point $(x,y)\in\R^2/\Z^2$.

If $\al/\beta$ is in $\Q$, then every leaf is compact. If $\al/\beta$ is not in $\Q$, then every leaf
is non-compact and dense in $M$.

Let $\CF$ be a foliation on the manifold $M$ and let $T_\CF\subset TM$ be the sub-bundle of all
vectors tangent to leaves. Since a sub-manifold is uniquely determined by its tangent bundle, a foliation $\CF$ is uniquely determined by its tangent bundle $T_\CF$.
Not every sub-bundle $T$ of $TM$ is tangent to a foliation. A sub-bundle that is tangent to a foliation is called
\emph{integrable}. There are other characterizations of integrability. For instance, a bundle
$T\subset TM$ is integrable if and only if for any two vector fields $X,Y$ with
$X,Y\in\Ga^\infty(T)$ it follows that $[X,Y]\in\Ga^\infty(T)$. Let
$$
\Omega_\CF^p\df \Ga^\infty(\wedge^p T_\CF^*).
$$
Using foliation coordinates one shows that the exterior differential of the de Rham complex of $M$
induces a differential $d:\Omega_\CF^p\ra\Omega_\CF^{p+1}$. with $d^2=0$. The \emph{foliation
cohomology} is
$$
H^\bullet(\CF)\df \ker(d)/\im(d).
$$
For any $p$ the space $\Omega_\CF^p$ is a Fr\'echet space, but the differential $d$ does not in
general have closed image, which implies that the quotient topology may be non-Hausdorff. Thus it
seems natural to define the \emph{reduced foliation cohomology} as the corresponding Hausdorff
quotient, i.e.,
$$
\bar H^\bullet(\CF)\df \ker(d)/\overline{\im(d)}.
$$ 
Let $E$ be a vector bundle which has a flat connection along the leaves. Then we have a
differential on the $E$-valued differential forms and we can form the corresponding reduced
cohomology which we write as
$\bar H^\bullet(\CF\otimes E)$.

\section{Anosov flows}
A smooth flow $\phi_t$ on a compact manifold $M$ is called  an \emph{Anosov flow} \cite{Smale}, if
the tangent bundle $TM$ of $M$ splits as
$$
TM\= T_0\oplus T_s\oplus T_u,
$$
where $T_0$, the \emph{neutral bundle}, is of rank one and generated by the flow $\phi_t$. Note that
this implies that the flow has no fixed points.  Next, $T_s$, the \emph{stable bundle} consists of all
$v\in TM$ such that $\norm{\phi_{t,*}v}\ra\infty$ as $t\ra+\infty$ for any Riemannian metric on $M$.
Since $M$ is compact this property does not depend on the choice of the metric. Finally, $T_u$, the
\emph{unstable bundle} comprises all vectors $v\in TM$ such that $\norm{\phi_{t,*}v}\ra\infty$ as
$t\ra-\infty$.

It turns out that the bundles $T_u,T_s$ are integrable so there are corresponding foliations, the
unstable and the stable foliation. For instance, two points $m,n\in M$ lie in the same leaf of the
stable foliation if and only if $d(\phi_tm,\phi_tn)$ tends to zero as $t\ra +\infty$, where $d$ is the
distance function of any Riemannian metric on $M$.

Suppose that $\dim M>1$. Since the manifold is compact, both $T_u$ and $T_s$ have to be
nonzero, so the smallest dimension in which a nontrivial Anosov flow can exist is three. An
example is given as follows. Let $\Ga\subset G=\SL_2(\R)/\pm 1$ be a discrete cocompact
subgroup and set $M=\Ga\bs G$. Let $H=\left(\begin{array}{cc} 1 & {}\\ {} & -1\end{array}\right)$.
Then for every $t\in\R$ the matrix $\exp(tH)=\left(\begin{array}{cc} e^t & {}\\ {} &
e^{-t}\end{array}\right)$ can be considered an element of $G$. On $M$ we get the flow
$$
\phi_t(\Ga g)\=\Ga\, g\,\exp(tH).
$$
Then $\phi$ is Anosov. The stable leaf through $\Ga g$ is given by $l_g= \Ga g N$,
where $N=\left\{\left.\pm\left(\begin{array}{cc}1 & x\\ {} & 1\end{array}\right)\right| x\in\R\right\}$.

Another way to get Anosov flows is to suspend a Anosov diffeomorphism. A diffeomorphism
$f:M\ra M$ is called Anosov if the tangent bundle $TM$ decomposes into a sum $T_s\oplus T_u$
of a stable and an unstable bundle which are defined as in the flow case. The suspension of an
Anosov diffeomorphism is an Anosov flow.

Examples of Anosov diffeomorphisms are constructed as follows. Let $G$ be a semi-simple Lie
group and $\Ga\subset G$ a discrete cocompact subgroup. The quotient $M=\Ga\bs G$ is called a
\emph{nilmanifold}. Let $f:G\ra G$ be an automorphism with $f(\Ga)=\Ga$. Then $f$ induces a
diffeomorphism on $M$ denoted by the same letter. Suppose that the differential $f_*:T_eG\ra
T_eG$ at the unit element $e$ of $G$ is hyperbolic in the sense that for every eigenvalue $\la\in\C$
of $f_*$ we have $|\la|\ne 1$. Then the induced diffeomorphism is Anosov. Such a
diffeomorphism is called an \emph{algebraic Anosov diffeomorphism}. There is a conjecture
\cite{Margulis} that states that up to finite covering every Anosov diffeomorphism should be
topologically conjugate to an algebraic one. This means that up to finite covering  for every Anosov diffeomorphism $f$ on a smooth manifold $M$ there should be an algebraic Anosov diffeomorphism $g$ on some nil-manifold $\Ga\bs G$ and a homeomorphism $\ph\colon M\ra\Ga\bs G$ such that $f=\ph^{-1} g\ph$.

The following conjecture was, in a slightly different setting and formulation, first given by V.
Guillemin \cite{Guillemin} and later by S. Patterson \cite{Patterson}. For a flow, a \emph{closed orbit}
is considered to come with a multiplicity, so if you have a closed orbit $c$, then you can go through it
twice and get a different closed orbit $c^2$. So for every closed orbit $c$ there is an underlying
\emph{primitive} closed orbit $c_0$ such that $c$ is a power of $c_0$ but $c_0$ is not the power of
a shorter orbit.

\begin{conjecture}
Let $\phi_t$ be a Anosov flow with stable foliation $\CF_s$. Then, as an identity of distributions
on $(0,\infty)$ we have
$$
\sum_{p=0}^{\rank\CF_s}(-1)^p \tr(\phi_t^*\mid \bar H^p(\CF_s))\=\sum_{c\ \rm closed}
l(c_0)\frac{\delta(t-l(c))}{\det(1-\phi_{l(c)}^*\mid T_{s,x})},
$$
where on the right hand side $x$ is any point on the the orbit $c$.
\end{conjecture}

In the case of a flow which is suspended from an algebraic diffeomorphism, C. Deninger and the
author have proved this conjecture \cite{Deitmar-Deninger}.

\begin{theorem} (C. Deninger-AD)\\
The Guillemin-Patterson conjecture is true for flows suspended from algebraic Anosov
diffeomorphisms. More specifically, if $f:M\ra M$ is an algebraic Anosov diffeomorphism with
stable bundle $\CF_s$, then the reduced cohomology $\bar H^\bullet(\CF_s)$ is finite dimensional
and
$$
\sum_{p=0}^{\rank\CF_s} (-1)^p\,\tr(f^*\mid \bar H^p(\CF_s))\=\sum_{x=f(x)}\frac{\det(1-f_*\mid
T_{s,x})}{|\det(1-f_*\mid T_xM)|}.
$$
\end{theorem}

For the proof one shows that the foliation cohomology is canonically isomorphic to Lie algebra
cohomology with trivial coefficients. This is shown inductively using the Hochschild-Serre spectral
sequence interatedly.

\section{The Selberg zeta function}
The Selberg zeta function for a compact Riemannian surface $Y$ of genus $g\ge 2$ is defined by
$$
Z_Y(s)\df \prod_{c_0}\prod_{N\ge 0} \left( 1-e^{-(s+N)l(c_0)}\right),
$$
where the first product runs over all primitive closed geodesics in $Y$ which is equipped with the
hyperbolic metric. Selberg showed in \cite{Selberg} that $Z_Y$ extends to an entire function which
satisfies a generalized Riemann hypothesis insofar as all zeros are in $\R\cup(\frac 12 i\R)$. 

In \cite{Cartier-Voros} P. Cartier and A. Voros gave the following determinant expression.

\begin{theorem}
(Cartier-Voros)\\
We have
$$
Z_Y(\frac 12+s)\=\left( e^{s^2}\det\left((\Delta_d+\frac 14)^{\frac
12}+s\right)\right)^{2g-2}\,\det\left((\Delta -\frac 14)+s^2\right).
$$
Here $\Delta$ is the Laplace operator on functions of $Y$ and $\Delta_d$ is the Laplace operator
on the sphere $S^2$.
\end{theorem}

This theorem can be proved by means of the trace formula inserting test functions of the form $f(\Delta)$ where $\Delta$ is the Laplace operator and $f$ a sufficiently nice function on the spectrum of $\Delta$. So either powers of the resolvent kernel, or heat or
wave kernels will do. These methods can be generalized to locally symmetric spaces of rank one
\cite{BuOl, det}, but they will not neatly generalize to higher rank because any
functional calculus of the given sort can not distinguish contributions of split tori of the same
dimension which are not conjugate. Since the zeta functions which are attached to such tori show
different analytical behaviour, a separation indeed is necessary.

Let now $Y$ denote a compact Riemannian manifold of odd dimension. Let
$SY$ denote the
\emph{sphere bundle} of
$Y$, i.e.,
$$
SY\df \{ v\in TY\mid \norm{v} =1\}.
$$
On $SY$ there is a natural flow, the \emph{geodesic flow} of $Y$. It is defined as follows. Let $t>0$.
A point $p$ of $SY$ comprises a point in $Y$ plus a direction. If you walk along the unique
geodesic in that direction for the time $t$, you get a new point and a new direction, the one you
came along in. Thus you get a new point $\phi_tp$ in $SY$. It is clear that closed orbits of the
geodesic flow correspond to closed geodesics in $Y$. It turns out that $\phi_t$ is Anosov. Let
$\CF_s$ denote its stable foliation. define the Selberg zeta function in this setting as
$$
Z_Y(s)\df \prod_{c_0}\prod_{N\ge 0} \det(1-e^{-sl(c_0)}\phi_{l(c_0)}\mid S^N(T_{s,x})),
$$
where $x$ is any point on the primitive closed orbit $c_0$ of the geodesic flow $\phi_t$ and $S^N$
denotes the $N$th symmetric power.

\begin{theorem}
(Lefschetz determinant formula, AD)
$$
Z_Y(s)\=\prod_{q\ge 0}\det(F+s\mid \bar H^q(\CF_s))^{(-1)^{q+1}}.
$$
\end{theorem}

This theorem can be proved as follows. First one uses the Lefschetz trace formula for rank one
spaces \cite{Juhl} to see that the divisors on both sides agree. The actual existence of the
determinants follows from the work of G. Illies \cite{Illies}. This proves the identity up to a factor of
the form $e^{P(s)}$, where $P$ is a polynomial. Finally one uses the asymptotic of the regularized
determinants \cite{l2det} to conclude the proof of the theorem.

\section{The Lefschetz formula for higher rank}
Let $G$ be a connected semisimple Lie group with finite center. Let $X=G/K$ be the attached
globally symmetric space, where $K\subset G$ is a maximal compact subgroup. Let $Y=\Ga\bs
X=\Ga\bs G/K$, where $\Ga$ is a discrete, cocompact, torsion-free subgroup of $G$. Then $Y$ is a
locally symmetric space. Let $SX$ and $SY$ denote the sphere bundles, then $SY=\Ga\bs SX$.
We want to understand the $G$-orbit structure of $SX$. For this recall that $G$ acts transitively on
$X=G/K$. So we get
$$
G\bs SX\= K\bs S_{eK}X\= W\bs S(A),
$$
where $A$ is a maximal split torus in $G$, so $A\cong \R^r$ as a Lie group, and $S(A)$ is the
sphere of norm one elements in $A$. Finally, $W$ denotes the Weyl group $W=N(A)/Z(A)$, the
quotient of the normalizer of $A$ and the centralizer of $A$. Then $W$ is a finite reflection group
acting on $A$.

The set $G\bs SX$ can be identified with the set of norm-1 elements of a closed positive Weyl chamber,
and so, $G\bs SX$ has the orbifold structure of a polysimplex and the latter can be viewed as a
subset of $A$ in a natural way. Let $f$ be a facet of this polysimplex and let $A_f\subset A$ be the
subgroup generated by $f$. We say that $f$ is \emph{cuspidal} if $A_f$ is the split part of a
Cartan subgroup of $G$. As an example consider $G=\SL_3(\R)$. Then $G\bs SX$ has dimension
one, so is a closed interval and has three facets, the open one and the two endpoints. In this case
each facet is cuspidal. For $G=\SL_4(\R)$ the polysimplex is the two dimensional simplex, so it has
$7$ facets. With the exception of one of the vertices each facet is cuspidal. Generally, the open
facet always is cuspidal and the bigger the dimension of a facet, the more likely it will be cuspidal.

Fix a cuspidal facet $f$. Let $C_f$ denote the set of all closed geodesics in $Y$ that lift into $f$.
Every $c\in C_f$ gives a point $a_c$ in the positive Weyl chamber $A_f^+$ by taking the
corresponding point in $S(A_f^+)$ and multiplying it with the length of $c$.

The pullback $Gf$  of the facet $f$ is, as a $G$-set, isomorphic to $(G/K_f)\times e$, where $e$ is a
cell on which $G$ acts trivially and $K_f=Z(A)\cap K$. So $A_f$ acts on $G/K_f$ and on $\Ga \bs
G/K_f$ by right multiplication. This action is Anosov in the sense that the tangent bundle of
$G/K_f$ or $\Ga\bs G/K_f$ decomposes as $T_0\oplus T_n\oplus T_s\oplus T_u$, where $T_n$ is
and additional 
\emph{neutral bundle} on which $A_f$ preserves norms, $T_s$ is the stable bundle, which
comprises all vectors that tend to zero under the positive Weyl cone. The unstable bundle $T_u$
finally consists of all vectors which tend to zero under the negative Weyl cone $A_f^-=\{ a^{-1}\mid
a\in A_f^+\}$. 

\begin{theorem}
(Lefschetz trace formula, AD)\\
Let $\CF_s$ denote the stable foliation.
As a distribution on $A_f^+$ we have that
$$
\sum_{p,q\ge 0} (-1)^{p+q+\rank \CF_s}\, \tr\left( a\mid \bar H^q(\CF_s\otimes\wedge^p T_n)\right)
$$
equals
$$
\sum_{c\ \rm
closed/homotopy}
\frac{\la_c\,\chi(A_f\bs X_c)\,\delta(a-a_c)}{\det(1-\phi_{l(c)}\mid T_{s,x})},
$$
where $\la_c$ is the volume of the unique compact $A_f$-orbit that contains $c$. Further $X_c$ is
the union of all closed geodesics homotopic to $c$ and $\chi(A_f\bs X_c)$ is the Euler-characteristic
of the quotient of $X_c$ modulo the $A_f$-action.
\end{theorem}

The proof \cite{hr} requires a geometric construction of a test function to be plugged into the trace formula.
This test function is built in a way as to have non-trivial orbital integrals only on conjugates of a
prescribed Cartan subgroup $H$. It is defined by $f(ghg^{-1})=\eta(g)\ph(h)$, where $\eta$ is a
suitable function on the homogeneous space $G/H$. The spectral interpretation in terms of foliation
cohomology is a consequence of the Osborne conjecture \cite{HeSch}.

\section{The higher rank Selberg zeta function}
Suppose there is a cuspidal facet $f$ of dimension zero, i.e., $f$ is a vertex of the polysimplex $G\bs
SX$. Then $A_f$ is one dimensional. Let $F$ be the positive infinitesimal generator of norm $1$.
Then
$F$ can also be viewed as the infinitesimal generator of the geodesic flow, i.e.,
$F=\frac{d}{dt}|_{t=0}\phi_t$. Define the \emph{generalized Selberg zeta function (AD)} as
$$
Z_f(s)\df \prod_{\stackrel{c\in C_f} {\rm prime}} \prod_{N\ge 0}\, \det\left( 1- e^{-sl(c)}\phi_{l(c),x}\mid
S^N(T_s)\right)^{\chi(A_f\bs X_c)}.
$$

\begin{theorem} (AD)
The double product defining $Z_f(s)$ converges if $\Re(s)>>0$. The function $Z_f(s)$ extends to a
meromorphic function  on the plane. It has finitely many poles which are located at real numbers
and under a suitable normalization of the metric all poles and zeros lie in $\R\cup(\frac 12 +i\R)$. The
vanishing order at
$s=\la$ equals
$$
(-1)^{\rank\CF_s}\sum_{p,q\ge 0} (-1)^{p+q}\, \dim \bar H^q(\CF_s\otimes \wedge^pT_n)_{\la-{\rm
eigenspace\ of}\ F}
$$
\end{theorem}

The proof uses the Lefschetz formula, resp. an extension to tempered distributions of the latter. One
plugs in a test function that has the property that the geometric side of the Lefschetz formula equals
a high logarithmic derivative of the Selberg zeta function.

\section{The higher rank prime geodesic theorem}
We will first state the classical prime geodesic theorem. Let $Y$ denote a compact Riemannian
surface of genus $\ge 2$ equipped with the hyperbolic metric.

\begin{theorem}
(Prime Geodesic theorem)\\For $x>0$ let
$$
\pi(x)\df \# \{ c_0\ {\rm prime}\mid e^{l(c_0)}\le x\}.
$$
Then under a suitable scaling of the metric,
$$
\pi(x)\ \sim\ \frac{x}{\log x},
$$
as $x\ra \infty$.
\end{theorem}

To motivate the higher rank case we will rewrite this theorem in the Chebysheff form. Define the
Chebysheff function by
$$
\psi(x)\df \sum_{c: e^{l(c)}\le x} l(c_0).
$$
Here the sum runs over all closed geodesics with the given length restriction and $c_0$ is the
prime geodesic underlying $c$. The prime geodesic theorem is equivalent to saying that, as
$x\ra\infty$,
$$
\psi(x)\ \sim\ x.
$$ 

Now let $X=G/K$ be a globally symmetric space as before and let $Y=\Ga\bs X$ be a compact
quotient. Let $f$ be the open facet of $G\bs SX$. The closed geodesics $c\in C_f$ are also called
the \emph{regular} geodesics. Every $c\in C_f$ gives a point $a_c\in A^+$. On $A^+$ there are
canonical coordinates stemming from primitive roots. To $c\in C_f$ we can thus attach coordinates
$c_1,\dots,c_r>0$.

\begin{theorem}
(Higher Rank Prime Geodesic Theorem, AD)\\
For $x_1,\dots,x_r>0$ let
$$
\psi(x_1,\dots,x_r)\df \sum_{\stackrel{c\colon e^{c_j}\le x_j}{j=1,\dots,r}} \la_c,
$$
where $\la_c$ is the volume of the unique maximal flat containing $C$. Then, as
$x_1,\dots,x_r\ra\infty$,
$$
\psi(x_1,\dots,x_r)\ \sim\ x_1\cdots x_r.
$$
\end{theorem}

To prove this theorem \cite{prime-geodesic} one uses a generalization of the logarithmic derivative of the Selberg zeta
function in several variables. as well as methods from analytical number theory (Tauberian
Theorems) extended to several variables.

We will close this section with a number theoretical application of the higher rank prime geodesic
theorem. This requires the prime geodesic theorem for locally symmetric manifolds $Y=\Ga\bs
X=\Ga\bs G/K$ which are not compact but of finite volume. In that case one needs to employ the
Arthur-Selberg trace formula which is quite harder to handle than the classic Selberg trace formula.
So it is not surprising that results in this setting are a lot more fragmentary at the moment.

In \cite{Sarnak} P. Sarnak proved the prime geodesic theorem for the arithmetic group $\SL_2(\Z)\bs
\SL_2(\R)/{\rm SO(2)}$ and inferred the following theorem.

\begin{theorem}
(Sarnak 83)\\
Let
$$
\pi_{2,0}(x)\df \sum_{\CO\colon e^{R(\CO)}\le x} h(\CO),
$$
where the sum ranges over all orders $\CO$ in real quadratic number fields, $R(\CO)$ denotes
the regulator of the order $\CO$, and $h(\CO)$ the class number of $\CO$. Then, as $x\ra\infty$,
$$
\pi_{2,0}(x)\ \sim\ \frac x{\log x}.
$$
\end{theorem}

Together with W. Hoffmann the author was recently able to prove a similar result for $\SL_3(\Z)$.

\begin{theorem}
(W Hoffmann-AD, 02)
Let 
$$ 
\pi_{1,1}(x)\df \sum_{\CO\colon e^{3R(\CO)}\le x} h(\CO),
$$
where the sum runs over all orders $\CO$ in number fields $F$ which have one real and two
complex embeddings. Then, as $x\ra\infty$,
$$
\pi_{1,1}(x)\ \sim\ \frac x{\log x}.
$$
\end{theorem}

For the proof one employs the Arthur-Selberg trace formula. First one plugs in functions that vanish
on parabolically degenerate elements to obtain a simplified trace formula whose geometric side
comprises orbital integrals only. Next one constructs test functions that are products of virtual
characters and twisted heat kernels to single out the relevant classes. The Mellin transform of the
resulting contribution is an analytic function bearing similarity to the logarithmic derivative of the
Selberg zeta function. The continuous spectral contributions cannot be computed but for the
purpose of the theorem it suffices to give a good estimate.

{\small University of Exeter, Mathematics, Exeter EX4
4QE, England\\ a.h.j.deitmar@ex.ac.uk}

\end{document}